\documentclass[letter]{article}
\usepackage[english]{babel}
\usepackage[latin1]{inputenc}
\usepackage[T1]{fontenc}
\usepackage{url}
\usepackage{graphicx}
\usepackage{subfigure}
\usepackage{color}
\usepackage{amsmath,amsfonts,amssymb}
\usepackage{colortbl}                   
\usepackage{url}
\usepackage{booktabs}
\usepackage{epsfig}
\usepackage{pstricks}

\usepackage{pst-node}
\usepackage{ifthen}
\usepackage{algorithm}
\usepackage{algpseudocode}
\usepackage{setspace}
\usepackage{pstricks-add}
\usepackage[framed]{ntheorem}
\usepackage{framed}
\usepackage{appendix}

\theoremstyle{definition}

\shadecolor{black!10!white}
\theorembodyfont{\rm\normalfont}
\newshadedtheorem{thm}{Theorem}

\newshadedtheorem{lem}{Lemma}
\newshadedtheorem{proposition}{Proposition}

\theoremstyle{nonumberplain}
\newtheorem{proof}{Proof}

\newcommand{\mypar}[1]{\vspace{2em}\noindent {\bf #1.}}

\newcommand{\qed}{\hfill $\Box$}

\title{A Proof of Convergence For the Alternating Direction Method of Multipliers Applied to Polyhedral-Constrained Functions}

\author{Jo\~ao F.~C.~Mota, Jo\~ao M.~F.~Xavier, Pedro M.~Q.~Aguiar, and~Markus P\"uschel}

\begin{document}

	\maketitle
	
	\begin{abstract}
      We give a general proof of convergence for the Alternating Direction Method of Multipliers (ADMM). ADMM is an optimization algorithm that has recently become very popular due to its capabilities to solve large-scale and/or distributed problems. We prove that the sequence generated by ADMM converges to an optimal primal-dual optimal solution. We assume the functions~$f$ and~$g$, defining the cost~$f(x) + g(y)$, are real-valued, but constrained to lie on polyhedral sets~$X$ and~$Y$. Our proof is an extension of the proofs from~\cite{Bertsekas:Parallel,BoydADMM}.
  \end{abstract}

  \section{Introduction}
  
  The Alternating Direction Method of Multipliers (ADMM), first proposed in the seventies by~\cite{Glowinski75-ADMMFirst,Gabay76-ADMMFirst}, is a versatile algorithm that is well-suited to solve large-scale or distributed problems. Lately, ADMM has become very popular because it efficiently handles problems that cannot be solved by the conventional interior point methods. Although ADMM and some of its variations have been extensively studied (e.g., \cite{Eckstein92-DouglasRachfordSplittingMethod,Fukushima92-ADMMSeparableConvexProblems,Bertsekas:Parallel}), there are still some theoretical aspects left to explore. As an example, only very recently it has been proved that ADMM has an~$O(1/k)$ convergence rate~\cite{He11-OnTheConvergenceRateADM}. We refer to the very complete survey~\cite{BoydADMM} for recent applications, history, and extensions of ADMM.
  
  ADMM solves the following optimization problem:
  \begin{equation}\label{Eq:ADMMProb}
    \begin{array}{ll}
      \text{minimize} & f(x) + g(y) \\
      \text{subject to} & x \in X, \, y \in Y \\
                        & Ax + By = c\,,
    \end{array}
  \end{equation} 
  where the variable is~$(x,y)\in \mathbb{R}^{n_1} \times \mathbb{R}^{n_2}$, and~$f:\mathbb{R}^{n_1} \xrightarrow{} \mathbb{R}$, $g:\mathbb{R}^{n_2} \xrightarrow{} \mathbb{R}$ are given convex functions, $A \in \mathbb{R}^{m \times n_1}$, $B \in \mathbb{R}^{m \times n_2}$ are two given matrices, $X$ and~$Y$ are convex sets, and~$c$ is a constant. The augmented Lagrangian of~\eqref{Eq:ADMMProb} is
  $$
    L_{\rho}(x,y;\lambda) = f(x) + g(y) + \lambda^{\top}(Ax + By - c) + \frac{\rho}{2}\|Ax + By - c\|^2\,,
  $$
  where~$\rho>0$ is a predefined positive parameter, and~$\lambda \in \mathbb{R}^m$ is the dual variable associated to the constraint~$Ax + By = c$. ADMM solves~\eqref{Eq:ADMMProb} by concatenating the method of multipliers with one iteration of the nonlinear Gauss-Seidel algorithm~\cite{Bertsekas:Parallel}, i.e., it iterates the following equations on~$k$:
  \begin{equation}\label{Eq:ADMMAlg}
    \left\{
      \begin{array}{l}
        x^{k+1} \in \arg\min_{x \in X} L_{\rho}(x, y^{k};\lambda^k) \vspace{0,2cm}\\
        y^{k+1} \in \arg\min_{y \in Y} L_{\rho}(x^{k+1}, y; \lambda^k) \vspace{0,2cm}\\
        \lambda^{k+1} = \lambda^k + \rho(Ax^{k+1} + By^{k+1} - c)
      \end{array}\,.
    \right.
  \end{equation} 
  In words, the augmented Lagrangian~$L_{\rho}$ is first minimized with respect to (w.r.t.)~$x$, keeping~$y$ and~$\lambda$ fixed at~$y^k$ and~$\lambda^k$, respectively. Then, $L_\rho$ is minimized w.r.t.~$y$, but~$x$ is fixed at the new value~$x^{k+1}$ (and~$\lambda$ is fixed at~$\lambda^k$). Finally, the dual variable~$\lambda$ is updated in a gradient ascent way.
  
  There are many variations of ADMM and also many variations of the proof of its convergence. Here, we use the techniques of~\cite{Bertsekas:Parallel,BoydADMM} to prove a version of ADMM that applies to the case where~$X$ and~$Y$ are polyhedral, the functions~$f$ and~$g$ are real-valued and convex, and the matrices~$A$ and~$B$ have full column-rank. We prove that the sequence~$\{(x^k,y^k,\lambda^k)\}$ generated by~\eqref{Eq:ADMMAlg} has a single limit point and that this limit point is a primal-dual solution of~\eqref{Eq:ADMMProb} or, in other words, a saddle point of the augmented Lagrangian~$L_\rho$.

	\section{Proof of Convergence}

	We aim to prove the following theorem.	
	\begin{thm}[Convergence of ADMM]\label{thm:ConvergenceADMM}
	\label{thm:ADMMConvergence}
		Assume:
		\begin{enumerate}
			\item $f:\mathbb{R}^{n_1} \xrightarrow{} \mathbb{R}$ and~$g :\mathbb{R}^{n_2} \xrightarrow{} \mathbb{R}$ are convex functions over~$\mathbb{R}^{n_1}$ and~$\mathbb{R}^{n_2}$, respectively
			\item $X \subset \mathbb{R}^{n_1}$ and $Y \subset \mathbb{R}^{n_2}$ are polyhedral sets			
			\item Problem~\eqref{Eq:ADMMProb} is solvable (and denote its optimal objective by~$p^\star$)
      \item Matrices~$A$ and~$B$ have full column-rank
		\end{enumerate}
		Then,
		\begin{enumerate}
			\item $f(x^k) + g(y^k) \xrightarrow{} p^\star$ 
			\item $\{(x^k,y^k)\}$ has a single limit point~$(x^\star,y^\star)$; furthermore, $(x^\star,y^\star)$ solves~\eqref{Eq:ADMMProb}
			\item $\{\lambda^k\}$ has a unique limit point~$\lambda^\star$; furthermore, $\lambda^\star$ solves the dual problem of~\eqref{Eq:ADMMProb}:
			\begin{equation}\label{Eq:ADMMDualProb}
        \begin{array}{cl}
          \text{maximize} & F(\lambda) + G(\lambda) - \lambda^\top c \\
              \lambda &
        \end{array}\,,
			\end{equation}
			where~$F(\lambda) = \inf_{x \in X} (f(x) + \lambda^\top Ax)$ and~$G(\lambda) = \inf_{y \in Y} (g(y) + \lambda^\top By)$
		\end{enumerate}
	\end{thm}
	
	We present a proof for this theorem that is based on the proofs from~\cite{BoydADMM} and~\cite[Prop.4.2]{Bertsekas:Parallel}. While~\cite{BoydADMM} does not make assumption~$4$ and consequently does not prove claim~$2$, \cite{Bertsekas:Parallel} proves all our claims by assuming that~$B$ is the identity matrix. So, Theorem~\ref{thm:ADMMConvergence} generalizes the proof from~\cite{Bertsekas:Parallel} by considering more general matrices~$B$, and it also generalizes the proof of~\cite{BoydADMM} by introducing assumption~$4$ and proving claim~$2$. Note that, although our assumptions on~$f$ and~$g$ are more restrictive than the ones in~\cite{BoydADMM}, Theorem~\ref{thm:ConvergenceADMM} can be straightforwardly adapted to the~$f$  and~$g$ assumptions made in~\cite{BoydADMM}.
	
	Before proving Theorem~\ref{thm:ConvergenceADMM}, we need the following lemma.
	\begin{lem}\label{lem:MinFunctionMinGradient}
		Let~$\phi$ and~$\psi$ be two convex functions from~$\mathbb{R}^n$ to~$\mathbb{R}$, and let~$\psi$ be differentiable in~$\mathbb{R}^n$. Also, let~$X \subset \mathbb{R}^n$ be a closed convex set. Then,
		$$
			x^\star \in \arg\min_{x \in X} \,\,\phi(x) + \psi(x) 
			\qquad \Longleftrightarrow \qquad 
			x^\star \in \arg\min_{x \in X} \,\,\phi(x) + \nabla \psi (x^\star)^\top (x - x^\star)\,.
		$$
	\end{lem}
	\begin{proof}
		The optimality conditions~\cite[Prop.4.7.2]{Bertsekas:Convex} for the left-hand side are: there exists
		$
      d \in \partial \bigl(\phi (x^\star) + \psi (x^\star) \bigr)
    $
    such that~$d^\top (x - x^\star) \geq 0$ for all~$x \in X$. Since $\partial \bigl(\phi (x^\star) + \psi (x^\star) \bigr)= \partial(\phi (x^\star)) + \nabla \psi (x^\star)$, this coincides with the optimality conditions for the right-hand side.
		\qed
	\end{proof}

	Now we are in conditions to prove Theorem~\ref{thm:ConvergenceADMM}.
	\begin{proof}
    First, note that assumptions $1$-$3$ make proposition $5.2.1$ of~\cite{Bertsekas:Nonlinear} applicable, which says that strong duality holds: there is no duality gap for~\eqref{Eq:ADMMProb} and~\eqref{Eq:ADMMDualProb} and the dual problem~\eqref{Eq:ADMMDualProb} is solvable. The pair~$(x^\star,y^\star)$ will denote any primal solution (there exists at least one by assumption~$3$) and~$\lambda^\star$ will denote any dual solution (there exists at least one by strong duality). We will also use the following notation:
		$$
			p^{k} = f(x^{k}) + g(y^{k})\,,
			\qquad
			r^{k} = Ax^k + By^k - c\,.
		$$
		The proof consists of showing that the following inequalities hold (see the proofs below):
		\begin{align}
			&p^\star - p^{k+1} \leq {\lambda^{\star}}^\top r^{k+1}
			\label{Eq:ADMMProofIneq1}
			\\
		  &p^{k+1} - p^{\star} \leq -(\lambda^{k+1})^{\top} r^{k+1} - \rho (B(y^{k+1} - y^k))^{\top}(B(y^{k+1} - y^\star) - r^{k+1})
			\label{Eq:ADMMProofIneq2}
			\\
		  &V^{k+1} \leq V^k - \rho \|r^{k+1}\|^2 - \rho \|B(y^{k+1} - y^k)\|^2\,,
			\label{Eq:ADMMProofIneq3}
		\end{align} 
		where~$V^k$ is the Lyapunov function
		$$
			V^k := \frac{1}{\rho}\|\lambda^k - \lambda^\star\|^2 + \rho\|B(y^k - y^\star)\|^2\,.
		$$
		Once these inequalities are proven, point~$1$ of the theorem is proven due to the following. From~\eqref{Eq:ADMMProofIneq3}, $V^k \leq V^0$, meaning that~$\lambda^k$ and~$By^k$ are bounded. Furthermore, 
		$$
			\rho\sum_{k=0}^{\infty} \bigl(\|r^{k+1}\|^2 + \|B(y^{k+1} - y^k)\|^2 \bigr) \leq V^0\,.
		$$
		Since~$V^0$ is finite, $r^k \xrightarrow{} 0$ and $B(y^{k+1} - y^k) \xrightarrow{} 0$ as~$k\xrightarrow{} \infty$. These facts together with the fact that~$By^k$ is bounded, imply that~$p^k \xrightarrow{} p^\star$ (since the right-hand side of both~\eqref{Eq:ADMMProofIneq1} and~\eqref{Eq:ADMMProofIneq2} converge to zero). This proves claim~$1$. 
		
    We now prove that the sequence~$\{(x^k,y^k)\}$ is bounded, which implies that it has limit points. To see that, note that~$\{By^k\}$ is bounded and~$B$ has full column-rank (assumption~$4$), thus $\{y^k\}$ is also bounded. Also, $r^k = Ax^k + By^k - c \xrightarrow{} 0$, which implies that~$\{Ax^k\}$ is bounded; again, $A$ has full column-rank, meaning that~$\{x^k\}$ is bounded.
		
		Given that we know that the sequence~$\{(x^k,y^k)\}$ has limit points, we now observe that any of its limit points, say~$(\bar{x},\bar{y})$, is primal optimal. In fact, this limit point is optimal because~$p^k = f(x^k) + g(y^k) \xrightarrow{} p^\star$ and thus any subsequence of~$p^k$ also converges to~$p^\star$. Also, any subsequence of~$r^k = Ax^k + By^k - c$ converges to~$0$ and, together with the fact that the sets~$X$ and~$Y$ are closed, this implies that~$(\bar{x},\bar{y})$ is feasible in~\eqref{Eq:ADMMProb}.
		
		Note that, although any limit point of~$\{(x^k,y^k)\}$ is primal optimal, this sequence may not even converge. This is not the case because of assumption~$4$. We will prove this after proving inequalities~\eqref{Eq:ADMMProofIneq1}-\eqref{Eq:ADMMProofIneq3} and that any limit point of~$\{\lambda^k\}$ is dual optimal (note that~$\{\lambda^k\}$ is bounded and hence it has limit points).

		\mypar{Proof of~\eqref{Eq:ADMMProofIneq1}}
			We have seen that strong duality holds for the pair~\eqref{Eq:ADMMProb}, \eqref{Eq:ADMMDualProb}. Thus, $(x^{\star},y^{\star},\lambda^\star)$ satisfies the KKT conditions. In particular,
			$$
				(x^\star,y^\star) \in \arg\min_{x\in X,y\in Y} f(x) + g(y)+ {\lambda^\star}^\top (Ax + By - c)\,,
			$$
			which implies
			$$
				\underbrace{f(x^\star) + g(y^\star)}_{p^\star} + {\lambda^\star}^\top(\underbrace{Ax^\star + By^\star - c}_{=0}) \leq \underbrace{f(x^{k+1}) + g(y^{k+1})}_{p^{k+1}} + {\lambda^\star}^\top(\underbrace{Ax^{k+1} + By^{k+1} - c}_{r^{k+1}})\,,
			$$
			or
			$$
				p^\star - p^{k+1} \leq {\lambda^\star}^\top r^{k+1}\,.
			$$
			
		\mypar{Proof of~\eqref{Eq:ADMMProofIneq2}}
			To prove~\eqref{Eq:ADMMProofIneq2}, we start by working the optimization problems defining~$x^{k+1}$ and~$y^{k+1}$ in~\eqref{Eq:ADMMAlg}.
			\begin{align}
				x^{k+1} &\in \arg\min_{x\in X} L_\rho (x,y^k;\lambda^k) 
				\notag
				\\
				&=
				  \arg\min_{x\in X} f(x) + g(y^k) + {\lambda^k}^\top (Ax + By^k - c) + \frac{\rho}{2}\|Ax + By^k - c\|^2
				\notag
				\\
				\intertext{and using Lemma~\ref{lem:MinFunctionMinGradient}}
				&=
				  \arg\min_{x \in X} f(x) + (A^\top \lambda^k + \rho A^\top (Ax^{k+1} + By^k - c))^\top
				  (x - x^{k+1})
				\notag
				\\
				&=
					\arg\min_{x\in X} f(x) + (A^\top \lambda^k + \rho A^\top (Ax^{k+1} + By^k - c))^\top x
				\notag
				\\
				&=
				  \arg\min_{x\in X} f(x) + (\lambda^k + \rho (Ax^{k+1} + By^k - c))^\top Ax
				\notag
				\\
				\intertext{and since~$\lambda^{k+1} = \lambda^k + \rho r^{k+1}$}
				&=
				  \arg\min_{x\in X} f(x) + (\lambda^{k+1} - \rho B(y^{k+1} - y^k))^\top Ax
				\label{Eq:ADMMProof1}
			\end{align}
			Using the same reasoning,
			\begin{align}
				y^{k+1} &\in \arg\min_{y \in Y} L_\rho (x^{k+1},y;\lambda^k) 
				\notag
				\\
				&=
				  \arg\min_{y \in Y} f(x^{k+1}) + g(y) + {\lambda^k}^\top (Ax^{k+1} + By - c) + \frac{\rho}{2}\|Ax^{k+1} + By - c\|^2
				\notag
				\\
				&=
				  \arg\min_{y \in Y} g(y) + (B^\top \lambda^k + \rho B^\top (Ax^{k+1} + By^{k+1} - c))^\top(y - y^{k+1})
				\notag
				\\
				&=
				  \arg\min_{y \in Y} g(y) + (\lambda^k + \rho (Ax^{k+1} + By^{k+1} - c))^\top By
				\notag
				\\
				\intertext{and since~$\lambda^{k+1} = \lambda^k + \rho r^{k+1}$}
				&=
				  \arg\min_{y \in Y} g(y) + {\lambda^{k+1}}^\top By
				\label{Eq:ADMMProof2}
			\end{align}
			Now we apply~\eqref{Eq:ADMMProof1} and~\eqref{Eq:ADMMProof2} to~$x^{k+1}$ and~$x^\star$, and to~$y^{k+1}$ and~$y^\star$, respectively. We have
			\begin{align}
				f(x^{k+1}) + (\lambda^{k+1} - \rho B(y^{k+1} - y^{k}))^\top A x^{k+1}
				&\leq
				f(x^\star) + (\lambda^{k+1} - \rho B(y^{k+1} - y^{k}))^\top A x^\star
				\label{Eq:ADMMProof3}
				\\
				\intertext{and}
				g(y^{k+1}) + {\lambda^{k+1}}^\top By^{k+1}
				&\leq
				g(y^\star) + {\lambda^{k+1}}^\top By^\star\,.
				\label{Eq:ADMMProof4}
			\end{align}
			Summing up~\eqref{Eq:ADMMProof3} and~\eqref{Eq:ADMMProof4} we get
			\begin{align*}
				&
				  p^{k+1} + {\lambda^{k+1}}^\top(Ax^{k+1} + By^{k+1}) - \rho(B(y^{k+1} - y^k))^\top Ax^{k+1}
				\leq
				  p^\star + {\lambda^{k+1}}^\top (\underbrace{Ax^\star + By^\star}_{=c})- \rho(B(y^{k+1} - y^k))^\top Ax^\star
				\\
				&\Longleftrightarrow
				  p^{k+1} - p^\star \leq -{\lambda^{k+1}}^\top r^{k+1} - \rho (B(y^{k+1} - y^k))^\top A(x^\star - x^{k+1})
				\\
				\intertext{and since~$r^{k+1} = Ax^{k+1} + By^{k+1} - c$ and~$c = Ax^\star + By^\star$, we have~$r^{k+1} = A(x^{k+1} - x^\star) + B(y^{k+1} - y^\star)$, thus}
				&\Longleftrightarrow
					p^{k+1} - p^\star \leq -{\lambda^{k+1}}^\top r^{k+1} - \rho (B(y^{k+1} - y^k))^\top
					(B(y^{k+1} - y^\star) - r^{k+1})\,,
			\end{align*}
			which is inequality~\eqref{Eq:ADMMProofIneq2}.
			
		\mypar{Proof of~\eqref{Eq:ADMMProofIneq3}}
			We concatenate~\eqref{Eq:ADMMProofIneq1} and~\eqref{Eq:ADMMProofIneq2}:
			\begin{align}
				&
				  -{\lambda^\star}^\top r^{k+1} \leq -(\lambda^{k+1})^{\top} r^{k+1} - \rho (B(y^{k+1} - y^k))^{\top}(B(y^{k+1} - y^\star) - r^{k+1})
				\notag
			  \\
			  &\Longleftrightarrow
			    (\lambda^{k+1} - \lambda^{\star})^\top r^{k+1}
			    - \rho (B(y^{k+1} - y^k))^{\top}r^{k+1}
			    + \rho (B(y^{k+1} - y^k))^{\top}(B(y^{k+1} - y^\star)) \leq 0
			  \notag
			  \\
			  &\Longleftrightarrow
					2(\lambda^{k+1} - \lambda^{\star})^\top r^{k+1}
			    - 2\rho (B(y^{k+1} - y^k))^{\top}r^{k+1}
			    + 2\rho (B(y^{k+1} - y^k))^{\top}(B(y^{k+1} - y^\star)) \leq 0\,.
				\label{Eq:ADMMProof5}
			\end{align}
			All we have to do now is to manipulate~\eqref{Eq:ADMMProof5} in order to get~\eqref{Eq:ADMMProofIneq3}. Taking into account that~$\lambda^{k+1} = \lambda^{k} + \rho r^{k+1}$, the first term of~\eqref{Eq:ADMMProof5} becomes:
			\begin{align}
				  2(\lambda^{k+1} - \lambda^\star)^\top r^{k+1}
				&=
				  2(\lambda^k - \lambda^\star)^\top r^{k+1} + 2\rho\|r^{k+1}\|^2
				\notag
				\\
				&=
				  2(\lambda^k - \lambda^\star)^\top r^{k+1} + \rho\|r^{k+1}\|^2 + \rho\|r^{k+1}\|^2
				\notag
				\\
				\intertext{and replacing~$r^{k+1} = \frac{1}{\rho}(\lambda^{k+1} - \lambda^k),$}
				&=
				  \frac{2}{\rho}(\lambda^k - \lambda^\star)^\top (\lambda^{k+1} - \lambda^k) + \frac{1}{\rho}\|\lambda^{k+1} - \lambda^k\|^2 + \rho\|r^{k+1}\|^2
				\notag
				\\
				\intertext{and completing the square of the first two terms,}
				&=
				  \frac{1}{\rho}\|\lambda^k - \lambda^\star + \lambda^{k+1} - \lambda^k\|^2 - \frac{1}{\rho}\|\lambda^k - \lambda^\star\|^2 + \rho \|r^{k+1}\|^2
				\notag
				\\
				&=
				  \frac{1}{\rho}\Bigl(\|\lambda^{k+1} - \lambda^\star\|^2 - \|\lambda^k - \lambda^\star\|^2\Bigr) + \rho \|r^{k+1}\|^2\,.
				  \label{Eq:ADMMProof6}
			\end{align}
			Insert~\eqref{Eq:ADMMProof6} in~\eqref{Eq:ADMMProof5}:
			\begin{multline}\label{Eq:ADMMProof7}
				  \frac{1}{\rho}\Bigl(\|\lambda^{k+1} - \lambda^\star\|^2 - \|\lambda^k - \lambda^\star\|^2\Bigr) 
				+ \underbrace{
			    	\rho \|r^{k+1}\|^2 - 2\rho(B(y^{k+1} - y^k))^\top r^{k+1}
				  }_{E_1}
				\\+ 2\rho(B(y^{k+1} - y^k))^\top B(y^{k+1} - y^\star) 
				\leq 0\,.
			\end{multline} 
			Completing the square of~$E_1$, the last two terms of~\eqref{Eq:ADMMProof7} become
			\begin{align*}
				& 
				  \rho\|r^{k+1} - B(y^{k+1} - y^k)\|^2
				  -\underbrace{
						    \rho\|B(y^{k+1} - y^k)\|^2 + 2\rho(B(y^{k+1} - y^k))^\top B(y^{k+1} - y^\star)
									    }_{E_2}
				\\
				\intertext{and completing the square of~$E_2$,}
				=&
				  \rho\|r^{k+1} - B(y^{k+1} - y^k)\|^2
				  -\rho\Bigl(
						      \|B(y^{k+1} - y^k) - B(y^{k+1} - y^\star)\|^2
						      -
						      \|B(y^{k+1}-y^\star)\|^2
						    \Bigr)
				\\
				=&
				  \rho\|r^{k+1} - B(y^{k+1} - y^k)\|^2
				  +\rho\Bigl(
						      \|B(y^{k+1}-y^\star)\|^2
						      -
						      \|B(y^k-y^\star)\|^2
						    \Bigr)
			\end{align*}
			Replacing this in~\eqref{Eq:ADMMProof7},
			\begin{align}
				&
				    \frac{1}{\rho}\Bigl(\|\lambda^{k+1} - \lambda^\star\|^2 - \|\lambda^k - \lambda^\star\|^2\Bigr) 
				  +
				    \rho\|r^{k+1} - B(y^{k+1} - y^k)\|^2
				  +
				    \rho\Bigl(
						      \|B(y^{k+1}-y^\star)\|^2
						      -
						      \|B(y^k-y^\star)\|^2
						    \Bigr)
				  \leq 0
				\notag
				\\
				&\Longleftrightarrow
				  (\underbrace{
							\frac{1}{\rho}\|\lambda^{k+1} - \lambda^\star\|^2
							+
							\rho\|B(y^{k+1} - y^\star)\|^2
											}_{=V^{k+1}}
					)
					-
					(\underbrace{
							\frac{1}{\rho}\|\lambda^{k} - \lambda^\star\|^2
							+
							\rho\|B(y^{k} - y^\star)\|^2
											}_{=V^{k}}
					)
					\leq 
					-\rho\|r^{k+1} - B(y^{k+1} - y^k)\|^2 
				\notag
				\\
				&\Longleftrightarrow
				  V^k - V^{k+1} \geq \rho\|r^{k+1} - B(y^{k+1} - y^k)\|^2
				\notag
				\\
				&\Longleftrightarrow
				  V^k - V^{k+1} \geq \rho\|r^{k+1}\|^2 + \rho\|B(y^{k+1} - y^k)\|^2 - 2\rho (B(y^{k+1} - y^k))^\top r^{k+1}\,.
				\label{Eq:ADMMProof8}
			\end{align} 
			What is left to prove is that~$-2\rho (B(y^{k+1} - y^k))^\top r^{k+1} \geq 0$. This is derived  from~\eqref{Eq:ADMMProof4}:
			\begin{align}
				g(y^{k+1}) + {\lambda^{k+1}}^\top By^{k+1} &\leq g(y^k) + {\lambda^{k+1}}^\top B y^k 
				\label{Eq:ADMMProof9}
				\\
				g(y^k) + {\lambda^k}^\top B y^k &\leq g(y^{k+1}) + {\lambda^k}^\top By^{k+1}
				\label{Eq:ADMMProof10}
			\end{align}
			Adding~\eqref{Eq:ADMMProof9} and~\eqref{Eq:ADMMProof10},
			\begin{align*}
					(\lambda^{k+1} - \lambda^k)^\top By^{k+1} \leq (\lambda^{k+1} - \lambda^k)^\top By^k
				&\Longleftrightarrow
				  (\lambda^{k+1} - \lambda^k)^\top B(y^{k+1} - y^k) \leq 0
				\\
				&\Longleftrightarrow
				  \rho{r^{k+1}}^\top B(y^{k+1} - y^k) \leq 0
 				\\
 				&\Longleftrightarrow
 				  -2\rho(B(y^{k+1} - y^k))^\top r^{k+1} \geq 0\,.
			\end{align*}
			Therefore, from~\eqref{Eq:ADMMProof8} we get~\eqref{Eq:ADMMProofIneq3}:
			$$
				V^k - V^{k+1} \geq \rho\|r^{k+1}\|^2 + \rho\|B(y^{k+1} - y^k)\|^2\,.
			$$
			
		\mypar{Proof that any limit point of~\boldmath{$\{\lambda^k\}$} is dual optimal} 
		We had seen that~$\{\lambda^k\}$ is bounded and thus it has limit points. Let~$\bar{\lambda}$ be any limit point and let~$\mathcal{K} \subset \mathbb{N}$ be the set of indices that yield that limit point, i.e., $\{\lambda^k\}_{k \in \mathcal{K}} \xrightarrow{} \bar{\lambda}$. We will show that~$\bar{\lambda}$ is dual optimal. 
		
		Let~$(\bar{x},\bar{y})$ be a limit point of~$\{(x^k,y^k)\}_{ k \in \mathcal{K}}$. We had seen that~$(\bar{x},\bar{y})$ existed (because $\{(x^k,y^k)\}$, and thus~$\{(x^k,y^k)\}_{k \in \mathcal{K}}$, is bounded) and that it was primal optimal. From now on, we will assume that~$\{(x^k,y^k)\}_{ k \in \mathcal{K}} \xrightarrow{} (\bar{x},\bar{y})$ (since~$(\bar{x},\bar{y})$ is just a limit point of~$\{(x^k,y^k)\}_{ k \in \mathcal{K}}$, if necessary, take a subsequence of~$\mathcal{K}$).
		
		Now, define
		$$
			\hat{\lambda}^{k} = \lambda^{k} - \rho B(y^{k} - y^{k-1})\,,
		$$
		and note that, because~$B(y^{k} - y^{k-1}) \xrightarrow{}0$ implies~$\{B(y^{k} - y^{k-1})\}_{k \in \mathcal{K}} \xrightarrow{} 0$, $\{\hat{\lambda}^k\}_{k \in \mathcal{K}}$ and $\{\lambda^k\}_{k \in \mathcal{K}}$ have the same limit point~$\bar{\lambda}$.
		From the definition of~$F(\lambda)$  and equation~\eqref{Eq:ADMMProof1}, 
		\begin{align}
			  F(\hat{\lambda}^{k}) &= \inf_{x \in X} f(x) + (\hat{\lambda}^{k})^\top (Ax)
			  \notag
			\\
			&= 
			  f(x^{k}) + (\hat{\lambda}^{k})^\top Ax^{k}
			\label{Eq:ADMMProofDualProb1}
			\\
			&\leq
			  f(x) + (\hat{\lambda}^{k})^\top Ax\,, \qquad \forall_{x \in X}\,.
			\label{Eq:ADMMProofDualProb2}
		\end{align}
		Similarly, from the definition of~$G(\lambda)$ and~\eqref{Eq:ADMMProof2},
		\begin{align}
			  G(\lambda^{k}) &= \inf_{y \in Y} g(y) + (\lambda^{k})^\top By
			  \notag
			\\
			&=
			  g(y^{k}) + (\lambda^{k})^\top By^{k}
			\label{Eq:ADMMProofDualProb3}
			\\
			&\leq
			  g(y) + (\lambda^{k})^\top By\,, \qquad \forall_{y \in Y}\,.
			\label{Eq:ADMMProofDualProb4}
		\end{align}
		Adding equations~\eqref{Eq:ADMMProofDualProb1} and~\eqref{Eq:ADMMProofDualProb3} and taking the limit~$k\xrightarrow{} +\infty$ ($k \in \mathcal{K}$) on both sides,
		\begin{equation}\label{Eq:ADMMProofDualProb5}
			  \lim_{\begin{subarray}{c} 
								k\xrightarrow{}+\infty \\
								k \in \mathcal{K}
						  \end{subarray}	
				} 
				(F(\hat{\lambda}^{k}) + G(\lambda^{k}))
			  =
			  f(\bar{x}) + g(\bar{y}) + \bar{\lambda}^\top (A\bar{x} + B\bar{y})
				=
				p^\star + \bar{\lambda}^\top c
				=
			  L(\lambda^\star) + \bar{\lambda}^\top c\,,
		\end{equation}
		where the second-to-last equality follows from optimality and primal feasibility of~$(\bar{x},\bar{y})$, and the last equality follows from strong duality. Adding equations~\eqref{Eq:ADMMProofDualProb2} and~\eqref{Eq:ADMMProofDualProb4} and taking the limit~$k\xrightarrow{} +\infty$ ($k \in \mathcal{K}$) on both sides,
		$$
				\lim_{\begin{subarray}{c} 
								k\xrightarrow{}+\infty \\
								k \in \mathcal{K}
						  \end{subarray}
				} (F(\hat{\lambda}^{k}) + G(\lambda^{k}))
			  \leq
			  f(x) + g(y) + \bar{\lambda}^\top (Ax + By)\,, \qquad \forall_{x \in X}\,\forall_{y \in Y}\,.
	  $$
		In particular, we can take the infimum on the right-hand side:
		\begin{equation}\label{Eq:ADMMProofDualProb6}
			 \lim_{\begin{subarray}{c} 
								k\xrightarrow{}+\infty \\
								k \in \mathcal{K}
						  \end{subarray}
			} (F(\hat{\lambda}^{k}) + G(\lambda^{k}))
			  \leq
		   F(\bar{\lambda}) + G(\bar{\lambda})\,.
		\end{equation}
		Inequality~\eqref{Eq:ADMMProofDualProb6} and equation~\eqref{Eq:ADMMProofDualProb5} yield
		$$
			F(\bar{\lambda}) + G(\bar{\lambda}) - \bar{\lambda}^\top c \geq L(\lambda^\star)\,,
		$$
		showing that~$\bar{\lambda}$ solves~\eqref{Eq:ADMMDualProb}, being dual optimal.
		
	\mypar{Proof of convergence of the sequence~\boldmath{$\{(x^k,y^k,\lambda^k)\}$} generated by ADMM}
    Let~$(\bar{x},\bar{y},\bar{\lambda})$ be any limit point of the sequence~$\{(x^k,y^k,\lambda^k)\}$ generated by ADMM. We have proved that~$(\bar{x}, \bar{y})$ is primal optimal and~$\bar{\lambda}$ is dual optimal. Now we prove that~$\{(x^k,y^k,\lambda^k)\}$ has a unique limit point, being convergent.
    
    Since~$(\bar{x},\bar{y},\bar{\lambda})$ is a primal-dual solution, the inequalities~\eqref{Eq:ADMMProofIneq1}-\eqref{Eq:ADMMProofIneq3} hold with~$(x^\star,y^\star,\lambda^\star)$ replaced by~$(\bar{x},\bar{y},\bar{\lambda})$. And, if we had seen that~$V^k$ in~\eqref{Eq:ADMMProofIneq3} was convergent (because it was bounded and non-increasing), now we see that its limit is~$0$ if we use~$(\bar{x},\bar{y},\bar{\lambda})$ as a solution. This implies that both~$\lambda^k \xrightarrow{} \bar{\lambda}$ and~$B(y^k - \bar{y})\xrightarrow{}0$. Since~$B$ has full column-rank, we also have~$y^k \xrightarrow{} \bar{y}$. Since~$r^k - c = A(x^k - \bar{x}) + B (y^k - \bar{y}) \xrightarrow{} 0$ and~$A$ has full column-rank, also~$x^k \xrightarrow{} \bar{x}$.
    
    This shows that all the sequences~$x^k$, $y^k$, and~$\lambda^k$, produced by ADMM, converge.

		\hfill
		\qed
	\end{proof}

\bibliographystyle{amsplain}
\bibliography{admm}

\providecommand{\bysame}{\leavevmode\hbox to3em{\hrulefill}\thinspace}
\providecommand{\MR}{\relax\ifhmode\unskip\space\fi MR }
\providecommand{\MRhref}[2]{%
  \href{http://www.ams.org/mathscinet-getitem?mr=#1}{#2}
}
\providecommand{\href}[2]{#2}
\begin{thebibliography}{1}

\bibitem{Bertsekas:Nonlinear}
D.~{Bertsekas}, \emph{Nonlinear programming}, 2 ed., Athena Scientific, 1999.

\bibitem{Bertsekas:Convex}
D.~{Bertsekas}, A.~{Nedi$\acute{c}$}, and A.~{Ozdaglar}, \emph{Convex analysis
  and optimization}, Athena Scientific, 2003.

\bibitem{Bertsekas:Parallel}
D.~{Bertsekas} and J.~{Tsitsiklis}, \emph{Parallel and distributed computation:
  Numerical methods}, Athena Scientific, 1997.

\bibitem{BoydADMM}
S.~{Boyd}, N.~{Parikh}, E.~{Chu}, B.~{Peleato}, and J.~{Eckstein},
  \emph{Distributed optimization and statistical learning via the alternating
  method of multipliers}, Foundations and Trends in Machine Learning \textbf{3}
  (2010), no.~1.

\bibitem{Eckstein92-DouglasRachfordSplittingMethod}
J.~{Eckstein} and D.~{Bertsekas}, \emph{On the douglas-rachford splitting
  method and the proximal point algorithm for maximal monotone operators},
  Mathematical Programming \textbf{55} (1992), 293--318.

\bibitem{Fukushima92-ADMMSeparableConvexProblems}
M.~{Fukushima}, \emph{Application of the alternating direction method of
  multipliers to separable convex programming problems}, Computational
  Optimization and Applications \textbf{1} (1992), 93--111.

\bibitem{Gabay76-ADMMFirst}
D.~{Gabay} and B.~{Mercier}, \emph{A dual algorithm for the solution of
  nonlinear variational problems via finite element approximations}, Computers
  and Mathematics with Applications \textbf{2} (1976), no.~1, 17--40.

\bibitem{Glowinski75-ADMMFirst}
R.~{Glowinski} and A.~{Marrocco}, \emph{Sur l'approximation, par éléments finis
  d'ordre un, et la résolution, par pénalisation-dualité, d'une classe de
  problèmes de dirichelet non linéaires}, Revue Française d'Automatique,
  Informatique, et Recherche Opérationelle \textbf{9} (1975), no.~2, 41--76.

\bibitem{He11-OnTheConvergenceRateADM}
B.~{He} and X.~{Yuan}, \emph{On the~$\textrm{O}(1/t)$ convergence rate of
  alternating direction method},
  \url{http://www.optimization-online.org/DB_HTML/2011/09/3157.html}, 2011.

\end{thebibliography}

\end{document}